\documentclass[12pt]{amsart}       
\makeatletter
\@namedef{subjclassname@2020}{\textup{2020} Mathematics Subject Classification}
\makeatother
\pdfoutput=1
\usepackage{subfig}
\usepackage{txfonts}
\usepackage [latin1]{inputenc}
\usepackage{amssymb}
\usepackage{eucal}
\usepackage{graphicx}
\usepackage{amsmath}
\usepackage{amscd}
\usepackage[all]{xy}           
\usepackage{tikz}
\usepackage{amsfonts,latexsym}
\usepackage{xspace}
\usepackage{epsfig}
\usepackage{float}
\usepackage{color}
\usepackage{fancybox}
\usepackage{colordvi}
\usepackage{multicol}
\usepackage{colordvi}
\usepackage{ifpdf}
\ifpdf
  \usepackage[colorlinks,final,backref=page,hyperindex]{hyperref}
\else
  \usepackage[colorlinks,final,backref=page,hyperindex,hypertex]{hyperref}
\fi
\usepackage[active]{srcltx} 

\usepackage{graphicx}
\usepackage{epstopdf}
\usepackage{epsfig}

\def\a{\alpha}
\def\b{\beta}

\def\c{\cdot}

\def\D{\Delta}

\def\o{\otimes}

\def\v{\varepsilon}

\newtheorem{theorem}{Theorem}[section]
\newtheorem{prop}[theorem]{Proposition}
\theoremstyle{definition}
\newtheorem{defn}[theorem]{Definition}
\newtheorem{lemma}[theorem]{Lemma}
\newtheorem{coro}[theorem]{Corollary}
\newtheorem{prop-def}{Proposition-Definition}[section]
\newtheorem{coro-def}{Corollary-Definition}[section]

\newtheorem{remark}[theorem]{Remark}

\newtheorem{exam}[theorem]{Example}

\newcommand{\nc}{\newcommand}
\nc{\tred}[1]{\textcolor{red}{#1}}
\nc{\tblue}[1]{\textcolor{blue}{#1}}
\nc{\tgreen}[1]{\textcolor{green}{#1}}
\nc{\tpurple}[1]{\textcolor{purple}{#1}}
\nc{\btred}[1]{\textcolor{red}{\bf #1}}
\nc{\btblue}[1]{\textcolor{blue}{\bf #1}}
\nc{\btgreen}[1]{\textcolor{green}{\bf #1}}
\nc{\btpurple}[1]{\textcolor{purple}{\bf #1}}
\nc{\NN}{{\mathbb N}}
\nc{\ncsha}{{\mbox{\cyr X}^{\mathrm NC}}} \nc{\ncshao}{{\mbox{\cyr
X}^{\mathrm NC}_0}}

\nc{\bfk}{\mathbf{k}}

\newcommand{\efootnote}[1]{}

\newcommand{\delete}[1]{}

\nc{\mlabel}[1]{\label{#1}}  
\nc{\mcite}[1]{\cite{#1}}  
\nc{\mref}[1]{\ref{#1}}  
\nc{\meqref}[1]{\eqref{#1}}  
\nc{\mbibitem}[1]{\bibitem{#1}} 

\delete{
\nc{\mlabel}[1]{\label{#1}  
{\hfill \hspace{1cm}{\small\tt{{\ }\hfill(#1)}}}}
\nc{\mcite}[1]{\cite{#1}{\small{\tt{{\ }(#1)}}}}  
\nc{\mref}[1]{\ref{#1}{{\tt{{\ }(#1)}}}}  
\nc{\meqref}[1]{\eqref{#1}{{\tt{{\ }(#1)}}}}  
\nc{\mbibitem}[1]{\bibitem[\bf #1]{#1}} 
}

\nc{\id}{\mathrm{id}}

\nc{\li}[1]{\tblue{#1}}
\nc{\lir}[1]{\tblue{\underline{Li:} #1}}

\nc{\mclip}[2]{{\centering
		\includegraphics[scale=#1]{#2}
}}

\nc{\ts}[1]{\textcolor{purple}{\underline{Tianshui:}#1}}
\nc{\hh}[1]{\textcolor{red}{\underline{Huihui:}#1}}

\begin{document}

\title{Rota-Baxter operators on cocommutative Hopf algebras and Hopf braces}
%
\author{Huihui Zheng}
\address{School of Mathematics and Information Science, Henan Normal University, Xinxiang 453007, China}
\email{huihuizhengmail@126.com}

\author{Li Guo}
\address{Rutgers University, Newark, NJ 07102, USA}
\email{liguo@rutgers.edu}

\author{Tianshui Ma}
\address{School of Mathematics and Information Science, Henan Normal University, Xinxiang 453007, China}
\email{matianshui@htu.edu.cn}

\author{Liangyun Zhang}
\address{Nanjing Agricultural University, Nanjing 210095, China}
\email{zlyun@njau.edu.cn}
\date{\today}
\begin{abstract}
This paper studies the relationship of Rota-Baxter operators on cocommutative Hopf algebras with Hopf braces and the Yang-Baxter equation, with emphasis on the embedding of cocommutative Hopf braces into Rota-Baxter Hopf algebras. Through Hopf braces, we establish a connection between relative Rota-Baxter operators on cocommutative Hopf algebras and bijective 1-cocycles. Finally, we introduce the notion of symmetric Hopf braces, and establish the relationship between symmetric Hopf braces and Rota-Baxter Hopf algebras.
\end{abstract}

\subjclass[2020]{16T05,
16T25, 
17B38, 
}
\keywords{Rota-Baxter operator, Hopf algebra, Hopf brace, Yang-Baxter equation, 1-cocycle}

\maketitle

\tableofcontents

\setcounter{section}{0}

\allowdisplaybreaks


\section{Introduction}
Several algebraic structures related to the quantum and classicla Yang-Baxter equations have appeared in recent years. 

In 2006, Rump introduced the notion of braces \mcite{Ru,Ru1} as a generalization of Jacobson radical rings, to study involutive set-theoretical solutions to the Yang-Baxter equation. It also has connections with regular subgroups and left orderable groups \mcite{BCV}, flat manifolds \mcite{Ru2}, Hopf-Galois extensions \mcite{Bac}. In 2014, this notion was reformulated by Ced\'{o}, Jespers and Okni\'{n}ski in \mcite{CJO}. In 2017, Guarnieri and Vendramin \mcite{GV} defined skew left braces which give non-involutive solutions to the Yang-Baxter equation~\mcite{BaR,Ya}.
Here a set $G$ with two binary operations $\cdot$ and $\circ$ is called a {\bf skew left brace} if $(G, \cdot)$ and $(G, \circ)$ are groups and the following identity holds.
$$
a \circ (bc) = (a \circ b) a^{-1} (a \circ c), \quad a, b, c\in G,
$$
where $a^{-1}$ denotes the inverse of $a$ in $(G,\c)$.

In 2017, Angiono, Galindo and Vendramin \mcite{AGV} introduced Hopf braces and Hopf cobraces. In \mcite{ZLMZ}, Hopf cobrace is further studied. A Hopf brace is a new algebraic structure related to the Yang-Baxter equation and generalizes Rump's braces.

On the other hand, Rota-Baxter operators on Lie algebras are the operator forms of the classical Yang-Baxter equation, known as the $\mathcal{O}$-operators~\mcite{BD,Ku,STS}; while Rota-Baxter 
operators on associative algebras are essential in the renormalization of quantum field theory, made apparent in the framework of Connes and Kreimer~\mcite{CK}. 

With motivation from Factorization Theorem of Semenov-Tian-Shansky for a Lie group that is fundamental in applications to integrable systems~\cite{FRS,RS1,RS2}, the notion of Rota-Baxter groups were recently introduced by Guo, Lang and Sheng in \mcite{GLS}. A {\bf Rota-Baxter group} is a group $G$ endowed with a map $B : G \rightarrow G$ satisfying the identity
$$
B(g)B(h) = B\big(gB(g)hB(g)^{-1}\big), \quad g, h \in G.
$$
Subsequently, the study of Rota-Baxter groups has been carried out further 
in a series of articles \mcite{BG,BG1,CS,CMS,JSZ}. Especially,  Bardakov and Gubarev \mcite{BG} proved that every Rota-Baxter group $(G, \cdot, B)$ gives rise to a skew left brace $(G, \cdot, \circ_B)$, where $x\circ_B y := xB(x)yB(x)^{-1} ~\text{for all}~x,y\in G$. Moreover, they demonstrated that every skew left brace can be embedded into a Rota-Baxter group.

In 2021, Goncharov introduced the notion of a Rota-Baxter operator on a cocommutative Hopf algebra in \mcite{Go}. A coalgebra map $B: H\rightarrow H$ is called a {\bf Rota-Baxter operator on a cocommutative Hopf algebra $H$} if
\begin{equation}
 B(x)B(y)=B\Big(x_{(1)}B\big(x_{(2)})yS(B(x_{(3)})\big)\Big), \quad x, y\in H.
\mlabel{eq:rboh}
\end{equation}
It generalizes the Rota-Baxter operators of weight 1 on groups and Lie algebras, in the sense that Rota-Baxter operator of weight 1 on a Lie algebra (resp. on a group) could be uniquely extended to a Rota-Baxter operator on the universal enveloping algebra (resp. on the group algebra). In~\mcite{LST},  Li, Sheng and Tang generalized this notion to relative Rota-Baxter operators on Hopf algebra and introduced the notion of post-Hopf algebras which generalizes post-Lie algebras and post-groups~\mcite{BGST}.
Various interconnections among relative Rota-Baxter Hopf algebras, post-Hopf algebras and Hopf braces have been discovered.

The purpose of present work is to further investigate such connections. The goal is to obtain finer results by specializing to Rota-Baxter operators from relative Rota-Baxter operators. Especially, we address the following embedding question:
since a skew brace can be embedded into a Rota-Baxter group~\mcite{BG} as noted above, whether a Hopf brace can be embedded into a Rota-Baxter Hopf algebra?

Overall, the paper is organized as follows. In Section \mref{sec:rboc}, we give compatibilities of some derived structures from Rota-Baxter operators on cocommutative Hopf algebras, and construct Rota-Baxter operators on cocommutative Hopf algebras via factorizations.

In Section \mref{sec:rboch}, we first give a direct connection from Rota-Baxter operators on cocommutative Hopf algebras to  cocommutative Hopf braces (Theorem~\mref{prop:rboHb}).
Conversely, we show that every cocommutative Hopf brace can be embedded into the induced Hopf brace of a Rota-Baxter Hopf algebra (Theorem~\mref{thm:embed}). There is a similar embedding of a cocommutative post-Hopf algebra into a Rota-Baxter Hopf algebra (Corollary~\mref{prop:rbopHA}).
In addition, we provide new ways to obtain solutions of the Yang-Baxter equation from Rota-Baxter operators on cocommutative Hopf algebras (Proposition~\mref{prop:rboYBE}). We finally construct a Hopf brace by a relative Rota-Baxter operator of Hopf algebra and give the connection between relative Rota-Baxter operators on Hopf algebras and bijective 1-cocycles (Proposition~\mref{prop:cocyclerRbo}).

In Section \mref{sec:shopf}, we introduce the notion of a symmetric Hopf brace and characterize it by identities and $H^{\rm op}$-module Hopf braces. Here $H^{\rm op}$ denotes the opposite Hopf algebra of a Hopf algebra $H$. Examples of symmetric Hopf braces are provided by Hopf algebras that allows factorizations, and their relationship with Rota-Baxter Hopf algebras is given.

We summarize the relations among the various structures in the following diagram.

\vspace{4mm}
\linethickness{0.4pt}
\ifx\plotpoint\undefined\newsavebox{\plotpoint}\fi 
{\small
\begin{picture}(200,132)(0,0)
\put(10.5,78.5){Hopf brace $(H,\cdot,\circ)$}
\put(280,78.25){Rota-Baxter Hopf algebra}
\put(295,0){bijective 1-cocycle}
\put(10,0){relative Rota-Baxter operator}
\put(278,125){symmetric Hopf brace}
\put(10.5,124){$H^{\rm op}$-module Hopf brace}
\put(120.75,127.25){\vector(1,0){152}}
\put(169,129){Prop.  \mref{prop:smHb}}
\put(276.75,82.75){\vector(-1,0){179}}
\put(170,71){Thm.  \mref{prop:rboHb}}

\put(75.75,75.75){\vector(0,-1){68}}
\put(73.5,9){\vector(0,1){68}}
\put(23,38.25){Prop.  \mref{prop:HbrRbo}}

\put(342.75,74){\vector(0,-1){68}}
\put(340.5,9){\vector(0,1){68}}
\put(345.25,38.25){Coro. \mref{coro:rborRbo}}
\put(282.25,38.25){Remark \mref{remark:cocyclerbo}}

\put(98,80.25){\vector(1,0){179}}
\put(170,87){Thm.  \mref{thm:embed}}

\put(140.2,3.2){\vector(1,0){154}}
\put(293.75,0.75){\vector(-1,0){152}}
\put(163,6){Prop.  \mref{prop:rRbococycle} and \mref{prop:cocyclerRbo}}

\put(75.75,87.75){\vector(0,1){32.75}}
\put(32,101){Defn. \mref{defn:moduleHb}}
\put(170,105){solutions of YBE}
\put(314.25,87){\vector(-3,1){60}}
\put(341.75,87.75){\vector(0,1){32.75}}
\put(345,101){Prop.  \mref{prop:sHbrbo}}
\put(268,103){Prop.  \mref{prop:rboYBE}}
\end{picture}
}
\smallskip

\vspace{4mm}

{\bf Notations.} Throughout the paper, let $\bfk$ be a fixed field. Unless otherwise specified, linearity, modules and tensor products are all taken over $\bfk$.
For brevity, we write a comultiplication $\D(c)$ as $c_{(1)}\o c_{(2)}$ without the summation sign. All Hopf algebras in this paper are cocommutative.

\section{Rota-Baxter operators on cocommutative Hopf algebra and their constructions}
\mlabel{sec:rboc}

In this section, we give some results of Rota-Baxter operators on cocommutative Hopf algebras and use factorizations of Hopf algebras to construct Rota-Baxter operators on Hopf algebras.

\subsection{Compatibility of derived Rota-Baxter Hopf algebras}\
\mlabel{ss:constr}

We first recall two fundamental processes with which a Rota-Baxter operator on a cocommutative Hopf algebra can derive another Rota-Baxter operator.

\begin{lemma}
\mcite{Go} Let $B$ be a Rota-Baxter operator on a cocommutative Hopf algebra $H$.
\begin{enumerate}
\item If $\varphi:H\rightarrow H$ is a bialgebra automorphism, then $B^{(\varphi)}:=\varphi \circ B\circ\varphi^{-1}$ is also a Rota-Baxter operator of $H$. \mlabel{it:2.1a}
\item Define $\widetilde{B}:H\rightarrow H$ by
\begin{equation}
	\widetilde{B}(x):=S(x_{(1)})B(S(x_{(2)})).
\mlabel{eq:trb}
\end{equation}
Then $\widetilde{B}$ is also a Rota-Baxter operator on $H$. \mlabel{it:2.1b}
\end{enumerate}
\mlabel{lemma:rboconstruct1}
\end{lemma}

Note that every bialgebra homomorphism from a Hopf algebra to another Hopf algebra is actually a Hopf algebra homomorphism. Thus, every bialgebra automorphism from a Hopf algebra to another Hopf algebra is already a Hopf algebra automorphism.

We now show the commutativity of the above two processes of deriving Rota-Baxter operators from a given one.

\begin{lemma}
Let $H$ be a cocommutative Hopf algebra, $\varphi$  a bialgebra automorphism of $H$, and $B$ a Rota-Baxter operator on $H$. Then $(\widetilde{B})^{(\varphi)} = \widetilde{B^{(\varphi)}}$.
\mlabel{lemma:rboproperty1}
\end{lemma}
\begin{proof}
For all $g\in H$, we have
\begin{eqnarray*}
(\widetilde{B})^{(\varphi)}(g)&=&\varphi(\widetilde{B}(\varphi^{-1}(g)))\\
&=&\varphi(S(\varphi^{-1}(g)_{(1)})B(S(\varphi^{-1}(g)_{(2)})))\\
&=&\varphi(S(\varphi^{-1}(g_{(1)}))BS\varphi^{-1}(g_{(2)}))\\
&=&S(g_{(1)})\varphi B S \varphi^{-1}(g_{(2)})\\
\hspace{6cm} &=&\widetilde{B^{(\varphi)}}(g). \hspace{7cm} \qedhere
\end{eqnarray*}
 \end{proof}

A new Rota-Baxter Hopf algebra can be induced from a Rota-Baxter cocommutative Hopf algebra via the following procedure which, when repeated, induces a sequence of Rota-Baxter Hopf algebras.

Let $(H,B)$ and $(H',B')$ be Rota-Baxter Hopf algebras. A {\bf homomorphism $g:(H,B)\rightarrow(H',B')$} is a Hopf algebra homomorphism such that $g\circ B=B'\circ g$. It is obvious that Rota-Baxter Hopf algebras and homomorphisms among them form a category, which is denoted by $\mathbf{rbHA}$.

\begin{lemma}
\mcite{Go} Let $(H,B)$ be a Rota-Baxter cocommutative Hopf algebra.
\begin{enumerate}
\item  Define
\begin{equation}
g\circ_B h:=g_{(1)}B(g_{(2)})hS(B(g_{(3)})), \quad T(g):=S(B(g_{(1)}))S(g_{(2)})B(g_{(3)}), \quad g,h\in H.
\mlabel{eq:desc}
\end{equation}
Then $H_B:=(H,\circ_B,\Delta,T)$ is a cocommutative Hopf algebra, called the {\bf descendent Hopf algebra} of $H$.
\mlabel{it:2.3a}
\item
\mlabel{lemma:rboproperty3}
Let $(H,B)$ be a Rota-Baxter cocommutative Hopf algebra. Then, for $T$ defined in Eq.~\meqref{eq:desc}, we have
\begin{equation}
	\mlabel{eq:rbopp}
	B(x_{(1)})B(T(x_{(2)}))=\v(x)1, \quad x\in H.
\end{equation}
\item The operator $B$ is also a Rota-Baxter operator on the cocommutative Hopf algebra $H_B$.
\mlabel{it:2.3b}
 \item The map $B:H_B\rightarrow H$ is a homomorphism of Rota-Baxter Hopf algebras.
     \mlabel{it:2.3c}
\end{enumerate}
\mlabel{lemma:desHA}
\end{lemma}

We next show that the derived Rota-Baxter operators from $B$ in Lemma~\mref{lemma:rboproperty1} have isomorphic descendent Hopf algebras.
We next show that the descendent Hopf algebras induced by $B$, $\widetilde{B}$ and $B^{(\varphi)}$ are isomorphic.

\begin{lemma}
Let $(H,B)$ be a Rota-Baxter cocommutative Hopf algebra and $\varphi:H\rightarrow H$ be a bialgebra automorphism. Then, $H_B\cong H_{\widetilde{B}}\cong H_{B^{(\varphi)}}$ as Hopf algebras.
\mlabel{lemma:rboproperty2}
\end{lemma}

\begin{proof}
Lemma \mref{lemma:rboconstruct1} shows that $\widetilde{B}$ and $B^{(\varphi)}$ are Rota-Baxter operators on $H$. So, by Lemma \mref{lemma:desHA}, $H_{\widetilde{B}}$ and $H_{B^{(\varphi)}}$ are Hopf algebras.

Note that the antipode $S:H\rightarrow H$ is bijective since $S^2=\id$. Moreover, for all $g,h\in H$,
\begin{eqnarray*}
S(g\circ_B h)&=&S(g_{(1)}B(g_{(2)})hS(B(g_{(3)})))\\
&=&B(g_{(1)})S(h)S(B(g_{(2)}))S(g_{(3)})\\
&=&S(g_{(1)})g_{(2)}B(g_{(3)})S(h)S(g_{(4)}B(g_{(5)}))
\\
&=&S(g_{(1)})\widetilde{B}(S(g_{(2)}))S(h)S(\widetilde{B}(S(g_{(3)}))) \qquad (\text{by the definition of } \widetilde{B} \text{ and }  S^2=\id)\\
&=&S(g)\circ_{\widetilde{B}} S(h).
\end{eqnarray*}
So $S:H_B\rightarrow H_{\widetilde{B}}$ is an algebra homomorphism.

Since $H$ is cocommutative, $S$ is a coalgebra map. Then $S:H_B\rightarrow H_{\widetilde{B}}$ is also a bialgebra map. Hence $S$ is an  isomorphism from the Hopf algebra $H_B$ to the Hopf algebra $H_{\widetilde{B}}$.

Similarly,
\begin{eqnarray*}
\varphi(g)\circ_{B^{\varphi}} \varphi(h)&=&\varphi(g_{(1)})B^{\varphi}(\varphi(g_{(2)}))\varphi(h)S(B^{\varphi}(\varphi(g_{(3)})))\\
&=&\varphi(g_{(1)})\varphi (B(g_{(2)}))\varphi(h)S(\varphi(B(g_{(3)})))\\
&=&\varphi(g_{(1)}B(g_{(2)})hS(B(g_{(3)})))\\
&=&\varphi(g\circ_B h).
\end{eqnarray*}
 Hence $\varphi$ is an isomorphism from the Hopf algebra $H_B$ to the Hopf algebra $H_{B^{(\varphi)}}$.
\end{proof}

\subsection{Rota-Baxter operators from factorizations of cocommutative Hopf algebras}\
\mlabel{ss:factor}

A common construction of Rota-Baxter operators on associative algebras and Lie algebras is by linear factorizations of these algebras. Such constructions extend to the level of Rota-Baxter operators on groups~\mcite{GLS}. We now construct Rota-Baxter operators on cocommutative Hopf algebras via factorizations of Hopf algebras.

\begin{theorem}
Let $G$ be a cocommutative Hopf algebra with a factorization $G=HLM$ into Hopf subalgebras $H, L$ and $M$, such that $HLM$ is isomorphic to $H\o L\o M$ as vector space. Let $C$ be a Rota-Baxter operator on $L$, and $\ell h=h\ell$, $mC(\ell)=C(\ell)m$, for all $h\in H, \ell\in L, m\in M$. Define a linear operator $B$ on $G$ by
\begin{equation}
B(h\ell m):=\v(h)C(\ell)S(m), \quad h\in H, \ell\in L, m\in M.
\mlabel{eq:prodop}
\end{equation}
\begin{enumerate}
\item
$B$ is a Rota-Baxter operator on $G$.
\mlabel{it:rbofact1}
\item Let $\circ_B$ be the product defined in Eq.~\meqref{eq:desc} and $L_C=(L,\circ_C)$ be the descendent Hopf algebra of the Rota-Baxter Hopf algebra $(L,C)$.
Then $G_B\cong H\o L_C\o M^{\rm op}$ as Hopf algebras.
\mlabel{it:rbofact2}
\end{enumerate}
\mlabel{thm:rbofactor}
\end{theorem}

\begin{proof}
\meqref{it:rbofact1}
We first check that $B$ satisfies the Rota-Baxter equation \eqref{eq:rboh}. For all $h,h'\in H, \ell,\ell'\in L, m,m'\in M$, since
\begin{eqnarray*}
mS(C(\ell))&=&m_{(1)}S(C(\ell_{(1)}))S(m_{(2)})C(\ell_{(2)})S(C(\ell_{(3)}))m_{(3)}\\
&=&m_{(1)}S(C(\ell_{(1)}))C(\ell_{(2)})S(m_{(2)})S(C(\ell_{(3)}))m_{(3)}\\
&=&S(C(\ell))m,
\end{eqnarray*}
we have
\begin{eqnarray*}
&&B((h\ell m)_{(1)}B((h\ell m)_{(2)})h'\ell' m'S(B((h\ell m)_{(3)})))\\
&=&B(h_{(1)}\ell_{(1)}m_{(1)}B(h_{(2)}\ell_{(2)}m_{(2)})h'\ell'm'S(B(h_{(3)}\ell_{(3)}m_{(3)})))\\
&=&B(h\ell_{(1)}m_{(1)}C(\ell_{(2)})S(m_{(2)})h'\ell'm'S(C(\ell_{(3)})S(m_{(3)})))\\
&=&B(h\ell_{(1)}C(\ell_{(2)})m_1S(m_{(2)})h'\ell'm'S(C(\ell_{(3)})S(m_{(3)})))\\
&=&B(h\ell_{(1)}C(\ell_{(2)})h'\ell'm'S(C(\ell_{(3)})S(m)))\\
&=&B(h\ell_{(1)}C(\ell_{(2)})h'\ell'm'S(C(\ell_{(3)})S(m)))\\
&=&B(hh'\ell_{(1)}C(\ell_{(2)})\ell'm'mS(C(\ell_{(3)})))\\
&=&B(hh'\ell_{(1)}C(\ell_{(2)})\ell'S(C(\ell_{(3)}))m'm)\\
&=&\v(hh')C(\ell_{(1)}C(\ell_{(2)})\ell'S(C(\ell_{(3)})))S(m)S(m')\\
&=&\v(hh')C(\ell)C(\ell')S(m)S(m')\\
&=&\v(hh')C(\ell)S(m)C(\ell')S(m')\\
&=&B(h\ell m)B(h'\ell'm').
\end{eqnarray*}
This proves Eq.\eqref{eq:rboh}.
Furthermore, since $C$ and $S$ are coalgebra homomorphisms, it follows that $B$ is also a coalgebra homomorphism. Hence, $B:G\rightarrow G$ is a Rota-Baxter operator on $H$.

\smallskip

\noindent
\meqref{it:rbofact2}
Indeed, for all $h,h'\in H, \ell,\ell'\in L, m,m'\in M$, we have
$$
\hspace{2cm} (h\ell m)\circ_B(h'\ell'm')=hh'\ell_{(1)}C(\ell_{(2)})\ell'S(C(\ell_{(3)}))m'm
=hh'(\ell\circ_C \ell')m'm. \hspace{2cm}\qedhere
$$
\end{proof}

\begin{prop}
\mlabel{prop:rboconstruct3}
Let $H$ and $K$ be cocommutative Hopf algebras. Assume that $H$ is a left $K$-module bialgebra with the module action $\triangleright$, and $G=H\#K$ with the smash product
$$
(h\#k)(h'\#k'): = h(k_{(1)}\triangleright h')\#k_{(2)}k',
\quad h,h'\in H,k,k'\in K.
$$
Let $C$ be a Rota-Baxter operator on $K$ and let $K_C$ be the descendent Hopf algebra.
\begin{enumerate}
\item The map
$$B: G\rightarrow G, \quad B(h\#k)=\v(h)C(k),$$
is a Rota-Baxter operator.
\mlabel{it:constr1}
\item
We have $G_B\cong H\# K_C$.
\mlabel{it:constr2}
\end{enumerate}
\end{prop}

\begin{proof}
\meqref{it:constr1}
In fact, for any $h,h'\in H, k,k'\in K$, we have
\begin{eqnarray*}
B((h\#k)_{(1)}B(h\#k)_{(2)}(h'\#k')S(B(h\#k)_{(3)}))&=&B((h_{(1)}\#k_{(1)})B(h_{(2)}\#k_{(2)})(h'\#k')S(B(h_{(3)}\#k_{(3)})))\\
&=&B((h\#k_{(1)}C(k_{(2)}))(h'\#k')S(C(k_{(3)})))\\
&=&B(h((k_{(1)}C(k_{(3)}))\triangleright h')\#k_{(2)}C(k_{(4)})k'S(C(k_{(5)})))\\
&=&\v(hh')C(k_{(1)}C(k_{(2)})k'S(C(k_{(3)})))\\
&=&\v(hh')C(k)C(k')\\
&=&B(h\#k)B(h'\#k').
\end{eqnarray*}
Thus $B(hk)B(h'k')=B((hk)_{(1)}B(hk)_{(2)}h'k'S(B(hk)_{(3)}))$, that is, $B$ is a Rota-Baxter operator on $G$.
\smallskip

\noindent

\meqref{it:constr2}
We first show that $H$ is a left $K_{C}$-module bialgebra with the module action
$$\unrhd: k\unrhd h:=(k_{(1)}C(k_{(2)}))\triangleright h, \quad k\in K, h\in H.$$
In fact, since $C$ is a Rota-Baxter
operator on $K$, we have
\begin{eqnarray*}
(k\circ_{C}\ell)\unrhd h&=&((k\circ_{C}\ell)_{(1)}C((k\circ_{C}\ell)_{(2)}))\triangleright h\\
&=&(k_{(1)}C(k_{(2)})\ell_{(1)}S_K(C(k_{(3)}))C(k_{(4)})C(\ell_{(2)}))\triangleright h\\
&=&(k_{(1)}C(k_{(2)})\ell_{(1)}C(\ell_{(2)}))\triangleright h\\
&=&(k_{(1)}C(k_{(2)}))\triangleright ((\ell_{(1)}C(\ell_{(2)}))\triangleright h)\\
&=&k\unrhd (\ell\unrhd h),
\end{eqnarray*}
for all $k,\ell\in K$ and $h\in H$. Similarly, we can verify the other conditions.

We next prove $G_B\cong H\# K_C$. This is because for all $h,h'\in H, k,k'\in K$, we have
\begin{eqnarray*}
(h\#k)\circ_B(h'\#k')&=&(h_{(1)}\#k_{(1)})B(h_{(2)}\#k_{(2)})(h'\#k')S(B((h_{(3)}\#k_{(3)})))\\
&=&(h\#k_{(1)}C(k_{(2)}))(h'\#k'S(C(k_{(3)})))\\
&=&h(k_{(1)}C(k_{(2)})\triangleright h')\#k_{(3)}C(k_{(4)})k'S(C(k_{(5)}))\\
&=&h((k_{(1)}C(k_{(2)}))\triangleright h'))\#(k_{(3)}\circ_C k')\\
\hspace{5cm}&=&h(k_{(1)}\unrhd h')\#k_{(2)}\circ_C k'. \hspace{5.5cm}\qedhere
\end{eqnarray*}
\end{proof}

\section{Rota-Baxter operators, Hopf braces, Yang-Baxter equation and 1-cocycles}
\mlabel{sec:rboch}

In this section, we study the relationship of Rota-Baxter operators on cocommutative Hopf algebras with Hopf braces and solutions of the Yang-Baxter equation. We also discuss the relationship between the more general relative Rota-Baxter operators on cocommutative Hopf algebras and Hopf braces. In addition, we establish the connection between relative Rota-Baxter operators on cocommutative Hopf algebras and bijective 1-cocycles.

\subsection{Rota-Baxter operators, Hopf braces and Yang-Baxter equation}\

Here we show that a Rota-Baxter operator on cocommutative Hopf algebra gives rise to a cocommutative Hopf brace. Conversely, every cocommutative Hopf brace is naturally embedded into to a Rota-Baxter cocommutative Hopf algebra. Further Rota-Baxter operators give solutions of the Yang-Baxter equation.

We first obtain a Hopf brace via a Rota-Baxter operator on a cocommutative Hopf algebra.

\begin{defn}\mcite{AGV}
A {\bf Hopf brace} structure over a coalgebra $(H,\Delta,\varepsilon)$ consists of the following data:

\begin{enumerate}
\item a Hopf algebra structure $(H, \cdot, 1, \Delta, \varepsilon, S)$ ($H$ for short),

\item a Hopf algebra structure $(H, \circ, 1_{\circ}, \Delta, \varepsilon, T)$ ($H_{\circ}$ for short),
\item the compatibility condition
\begin{equation}
a\circ(bc)=(a_{(1)}\circ b)S(a_{(2)})(a_{(3)}\circ c), ~a,b,c\in H.
\mlabel{eq:brace}
\end{equation}
\end{enumerate}
\end{defn}

A homomorphism $f:(M,\c_M,\circ_M)\rightarrow(N,\c_N,\circ_N)$ of Hopf braces is a linear map such that $f:M\rightarrow N$ and $f:M_{\circ}\rightarrow N_{\circ}$ are Hopf algebra homomorphisms (see \mcite{AGV}). Let $\mathbf{HB}$ denote the category of Hopf braces, and let $\mathbf{cocHB}$ denote its subcategory of cocommutative Hopf braces.

\begin{theorem}
\mlabel{prop:rboHb}
Let $(H,\cdot,B)$ be a Rota-Baxter cocommutative Hopf algebra, and $\varphi:H\rightarrow H$ a bialgebra automorphism. Then the following conclusions hold.
\begin{enumerate}
\item
 $(H,\cdot,\circ_B)$ and $(H,\cdot^{\rm op},\circ_B)$ are cocommutative Hopf braces.
\mlabel{it:3.1a}
\item $(H,\c^{\rm op},\circ_B)\cong (H,\cdot,\circ_{\widetilde{B}})$ and $ (H,\c,\circ_B)\cong (H,\cdot,\circ_{B^{(\varphi)}})$ as Hopf braces.
\mlabel{it:3.1b}
\item
 \mlabel{it:3.1c}
A Rota-Baxter Hopf algebra homomorphism $f:(H,\cdot,B)\to (H',\cdot',B')$ induces a homomorphism $f:(H,\cdot,\circ_B)\to (H',\cdot',\circ_{B'})$ of Hopf braces, yielding a functor
$$
\mathbf{rbHA}\rightarrow \mathbf{cocHB}, \ F(H,B):=(H,\c,\circ_B), ~~~~F(f):=f,
$$
from the category of Rota-Baxter Hopf algebras to the category of cocommutative Hopf braces.
    \end{enumerate}
\end{theorem}

\begin{proof} \eqref{it:3.1a} By Lemma \mref{lemma:rboproperty2}, $(H,\cdot,\circ_B)$ is a Hopf algebra. For $g,h,\ell\in G$, we have

\begin{eqnarray*}
(g_{(1)}\circ_B h)S(g_{(2)})(g_{(3)}\circ_B \ell)&=&g_{(1)(1)}B(g_{(1)(2)})hS(B(g_{(1)(3)}))S(g_{(2)})g_{(3)(1)}B(g_{(3)(2)})\ell S(B(g_{(3)(3)}))\\
&=&g_{(1)}B(g_{(2)})h\ell S(B(g_{(3)}))\\
&=&g\circ_B(h\ell).
\end{eqnarray*}
Hence $(H,\cdot,\circ_B)$ is a cocommutative Hopf brace.

The assumption that the Hopf algebra  $(H,\c,\D,S)$ is cocommutative implies that $(H,\c^{\rm op},\D, S)$ is also a cocommutative Hopf algebra.
Furthermore, for $g,h,l\in G$, we have
\begin{eqnarray*}
(g_{(1)}\circ_B h)\c^{\rm op} S(g_{(2)})\c^{\rm op}(g_{(3)}\circ_B \ell)&=&g_{(3)(1)}B(g_{(3)(2)})\ell S(B(g_{(3)(3)}))S(g_{(2)})g_{(1)(1)}B(g_{(1)(2)})hS(B(g_{(1)(3)}))\\
&=&g_{(1)}B(g_{(2)})\ell hS(B(g_{(3)}))\\
&=&g\circ_B(h\c^{\rm op}\ell).
\end{eqnarray*}
Hence $(H,\c^{\rm op},\circ_B)$ is a cocommutative Hopf brace.
\smallskip

\noindent
\eqref{it:3.1b} Since $S(h\c^{\rm op} g)=S(gh)=S(h)S(g)$, we find that $S:(H,\c^{\rm op})\rightarrow (H,\c)$ is an algebra map.
By \meqref{it:3.1a} and Lemma \mref{lemma:rboproperty2}, we get $(G,\c^{\rm op},\circ_B)\cong (G,\cdot,\circ_{\widetilde{B}})$ as Hopf braces.

By \eqref{it:3.1a}, we conclude that $(H,\cdot,\circ_{B^{(\varphi)}})$ is a Hopf brace. Further by Lemma \mref{lemma:rboproperty2}, it is easy to see that $(H,\cdot,\circ_{B}) \cong (H,\cdot,\circ_{B^{(\varphi)}})$ as Hopf braces.

\smallskip
\noindent
\eqref{it:3.1c} By \eqref{it:3.1a}, $(H,\c,\circ_B)$ and $(H',\c,\circ_{B'})$ are Hopf braces. We next prove that $f$ is a homomorphism of Hopf braces from $(H,\c,\circ_B)$ to $(H',\c,\circ_{B'})$.

As a matter of fact, because $f:(H,B)\rightarrow(H',B')$ is a homomorphism of Rota-Baxter Hopf algebra, we have a Hopf algebra homomorphism $f:H\rightarrow H'$ such that $f\circ B=B'\circ f$. Moreover, for any $h,g\in H$, we obtain
\begin{eqnarray*}
f(h\circ_B g)&=&f(h_{(1)}B(h_{(2)})gS_H(B(h_{(3)})))\\
&=&f(h_{(1)})f(B(h_{(2)}))f(g)f(S_H(B(h_{(3)})))\\
&=&f(h_{(1)})B'(f(h_{(2)}))f(g)S_{H'}(B'(f(h_{(3)})))\\
&=&f(h)_{(1)}B'(f(h)_{(2)})f(g)S_{H'}(B'(f(h)_{(3)}))\\
&=&f(h)\circ_{B'}f(g).
 \end{eqnarray*}
 Hence $F$ is a functor.
\end{proof}

\begin{coro}
 Let $(H,B)$ be a Rota-Baxter cocommutative Hopf algebra. If Im$B\subset Z(H)$ (the center of $H$), then the Hopf brace $(H,\cdot,\circ_B)$ is the trivial one.
\mlabel{trivialHb}
\end{coro}

\begin{proof} Indeed, for all $g,h\in H$,
$g\circ_B h=g_{(1)}B(g_{(2)})hS(B(g_{(3)}))=gh$.
\end{proof}

We recall some properties of Hopf braces for later use.

\begin{lemma}
\mlabel{lemma:Hbproperty}
\mcite{AGV} Let $(H, \cdot,\circ)$ be a cocommutative Hopf brace.  Let $H_{\circ}$ denote the Hopf algebra $(H, \circ, 1_{\circ}, \Delta, \varepsilon, T)$.
\begin{enumerate}
\item   $H$ is a left $H_{\circ}$-module bialgebra with
\begin{equation}
a\rightharpoonup b=S(a_{(1)})(a_{(2)}\circ b),  \ a,b\in H.
\mlabel{eq:hbmodule}
\end{equation}
\item For all $a,b\in H$,
\begin{equation}
a\circ b=a_{(1)}(a_{(2)}\rightharpoonup b),
\mlabel{eq:hbm}
\end{equation}
\begin{equation}
\mlabel{eq:hbe}
ab=a_{(1)}\circ(T(a_{(2)})\rightharpoonup b).
 \end{equation}
 \end{enumerate}
\end{lemma}

We next use cocommutative Hopf braces to construct Rota-Baxter Hopf algebras. Let $(H,\cdot,\circ)$ be a cocommutative Hopf brace.
Consider the tensor coalgebra $H\o H$
with multiplication and antipode given by
\begin{eqnarray*}
	(x\o y)\ast(z\o t)&=&x_{(1)}\circ z\o y(x_{(2)}\rightharpoonup t),\\
	S'(x\o y)&=&T(x_{(1)})\o T(x_{(2)})\circ (x_{(3)}S(y)),
\end{eqnarray*}
where $S$ and $T$ are antipodes of the Hopf algebras $(H,\cdot)$ and $(H,\circ)$ respectively, and $\rightharpoonup$ is defined by Eq.~\eqref{eq:hbmodule}.

Consider the sub-Hopf algebras
$$G=\{(h_{(1)}\o h_{(2)})\mid h\in H\} \text{ and } L=\{(h\o 1)\mid h\in H\}$$
of $H\o H$, noticing that if $(H,\c,\circ)$ is a Hopf brace, then $1=1_{\circ}$ (see\mcite{AGV}).

\begin{theorem}
\mlabel{thm:HbrbHA}
\begin{enumerate}
\item
With the notions given above, there is the factorization of Hopf algebras
$H\otimes H=G*L$.
\mlabel{it:hbr1}
\item
\mlabel{it:3.4a}
The splitting operator
\begin{equation}
B: H\o H\to H\o H, \quad B(g\ast \ell):=\varepsilon(g)S'(\ell),\ g\in G, \ell\in L,
\mlabel{eq:split}
\end{equation}
is a Rota-Baxter operator and takes the form
\begin{equation}
\mlabel{eq:rbform}
B(x\o y)=T(x)\circ y\o 1, \quad x, y\in H.
\end{equation}
\item
\mlabel{it:3.4b}
For a homomorphism of cocommutative Hopf braces $f:(H,\c,\circ)\to (H',\c,\circ)$, the map $f\o f: H\o H \to H\o H$ is a homomorphism of Rota-Baxter Hopf algebras, yielding a functor
$$J:\mathbf{cocHB}\rightarrow \mathbf{rbHA}, \
J((H,\c,\circ))=(H\o H,*,B), ~~~~J(f)=f\o f.
$$
\end{enumerate}
\end{theorem}

\begin{proof}
\eqref{it:hbr1} The factorization follows since $x\o y=(y_{(1)}\o y_{(2)})\ast (T(y_{(3)})\circ x\o 1)$ for $x,y\in H$.
\smallskip

\noindent
\eqref{it:3.4a} By ~\cite[Prop.  2]{Go},  $B$ is a Rota-Baxter operator on $H\otimes H$. We verify Eq.~\meqref{eq:rbform} as follows.
\begin{eqnarray*}
B(x\o y)&=&B((y_{(1)}\o y_{(2)})\ast (T(y_{(3)})\circ x\o 1))\\
&=&\v(y_{(1)}\o y_{(2)})S'(T(y_{(3)})\circ x\o 1)\\
&=&S'(T(y)\circ x\o 1)\\
&=&T((T(y)\circ x)_{(1)})\o T((T(y)\circ x)_{(2)})\circ ((T(y)\circ x)_{(3)}S(1))\\
&=&T(T(y)_{(1)}\circ x_{(1)})\o T(T(y)_{(2)}\circ x_{(2)})\circ (T(y)_{(3)}\circ x_{(3)})\\
&=&T(x_{(1)})\circ y_{(1)}\o T(x_{(2)})\circ y_{(2)}\circ T(y_{(3)})\circ x_{(3)}\ \ (T^2=id)\\
&=&T(x)\circ y\o 1.
\end{eqnarray*}

\noindent
\eqref{it:3.4b}
We first prove that, for the homomorphism $f$ of Hopf braces, the map $f\o f$ is a homomorphism of Rota-Baxter Hopf algebras from $(H\o H,B)$ to $(H'\o H',B')$.

First,  $f\o f$ is an algebra homomorphism: for any $x,y,z,t\in H$,
\begin{eqnarray*}
(f\o f)((x\o y)*(z\o t))&=&(f\o f)(x_{(1)}\circ z\o y(x_{(2)}\rightharpoonup t))\\
&=&f(x_{(1)}\circ z)\o f(y(x_{(2)}\rightharpoonup t))\\
&=&f(x_{(1)})\circ f(z)\o f(y)f(x_{(2)}\rightharpoonup t)\\
&=&f(x_{(1)})\circ f(z)\o f(y)(f(x_{(2)})\rightharpoonup f(t))\\
&=&f(x)_{(1)}\circ f(z)\o f(y)(f(x)_{(2)}\rightharpoonup f(t))\\
&=&(f(x)\o f(y))*(f(z)\o f(t))\\
&=&(f\o f)(x\o y)*(f\o f)(z\o t).
\end{eqnarray*}
It is also easy to prove that $f\o f$ is a coalgebra homomorphism.

Finally, we prove that $(f\o f)\circ B=B'\circ (f\o f)$:
\begin{eqnarray*}
(f\o f)B(x\o y)&=&(f\o f)(T_H(x)\circ y\o 1)\\
&=&f(T_H(x)\circ y)\o f(1)\\
&=&f(T_H(x))\circ f(y)\o 1\\
&=&T_{H'}(f(x))\circ f(y)\o 1\\
&=&B'(f(x)\o f(y))\\
&=&B'(f\o f)(x\o y).
\end{eqnarray*}

In summary, $f\o f$ is a homomorphism of Rota-Baxter Hopf algebras. Then it is direct to check that $J$ a functor.
\end{proof}

\begin{remark}
Let $(A, \Delta, \circ, T)$ be a cocommutative Hopf algebra. Assume that $(A, A)$ is a matched pair of cocommutative Hopf algebras \mcite{Ka} with actions $\rightharpoonup$ and $\leftharpoonup$ such that
$$a \circ b=(a_{(1)}\rightharpoonup b_{(1)}) \circ (a_{(2)}\leftharpoonup b_{(2)}), \quad a, b \in A.$$
Then by \cite[Prop.  3.2]{AGV}, we find that $(A, \cdot, \circ)$ is a cocommutative Hopf brace with
$$a\cdot b=a_{(1)}\circ (T(a_{(2)})\rightharpoonup b), \quad S(a)=a_{(1)}\rightharpoonup T(a_{(2)}), \quad a,b\in A.$$
Then by Theorem \mref{thm:HbrbHA}, there is a Rota-Baxter cocommutative Hopf algebra $(A\o A,B)$, where $B(x\o y)=T(x)\circ y\o 1$, for all $x,y\in A$, and the tensor coalgebra $A\o A$ has the multiplication $\ast$ and antipode $S'$ defined by
$$
(x\o y)\ast(z\o t)=x_{(1)}\circ z\o y_{(1)} \circ (T(y_{(2)})\rightharpoonup ((x_{(2)}\rightharpoonup T(x_{(3)}))\circ (T(x_{(4)}\rightharpoonup T(x_{(5)}))\rightharpoonup (x_{(6)}\circ t)))),
$$
$$
S'(x\o y)=T(x_{(1)})\o (T(x_{(2)})\circ y_{(1)})\rightharpoonup T(y_{(2)}), \quad x, y, z, t\in A.
$$
\mlabel{remark:2.20}
\end{remark}

Theorem~\mref{thm:HbrbHA} yields the following PBW type theorem.
\begin{theorem}
\begin{enumerate}
\item Every cocommutative Hopf brace can be embedded into a Rota-Baxter cocommutative Hopf algebra, regarded as a Hopf brace by Theorem~\mref{prop:rboHb}.
\mlabel{it:3.5a}
\item There is a natural transformation $\mathcal{T}:I\rightarrow FJ$, where $I$ is the identity functor on the category  $\mathbf{cocHB}$ of cocommutative Hopf braces, and $F$ and $J$ are defined in Theorem \mref{prop:rboHb} and Theorem \mref{thm:HbrbHA}, respectively.
\mlabel{it:3.5b}
\end{enumerate}
\mlabel{thm:embed}
\end{theorem}

However, the functors $F$ and $J$ do not give an equivalence between the categories $\mathbf{cocHB}$ and $\mathbf{rbHA}$, and $H\o H$ is not an enveloping Rota-Baxter Hopf algebra of the Hopf brace $H$. It would be interesting to obtain an enveloping Rota-Baxter Hopf algebra of a Hopf brace and establish a PBW theorem.

\begin{proof}
 \eqref{it:3.5a} Let $(H,\cdot,\circ)$ be a cocommutative Hopf brace. Then $(H\o H,B)$ is a Rota-Baxter cocommutative Hopf algebra by Theorem \mref{thm:HbrbHA}. By Theorem \mref{prop:rboHb}, $(H\o H)_B=(H\o H, \ast, \circ_B)$ is a Hopf brace.

Embed $H$ into $H\o H$ by $\psi:H\to H\o H, g\mapsto 1\o g$. Then, $\psi$ is an isomorphism of the Hopf braces $H$ with Im$(\psi)$, where Im$(\psi)$ is considered as a sub-Hopf brace of $(H\o H)_B=(H\o H, \ast, \circ_B)$.

In fact, for any $g,h\in H$, we have
$$
\psi(g)\ast\psi(h)=(1\o g)\ast(1\o h)=1\o g(1\rightharpoonup h)=1\o gh=\psi(gh),
$$
and
\begin{eqnarray*}
\psi(g)\circ_B\psi(h)&=&\psi(g)_{(1)}\ast B(\psi(g)_{(2)})\ast\psi(h)\ast S'(B(\psi(g)_{(3)}))\\
&=&(1\o g)_{(1)}\ast B((1\o g)_{(2)})\ast(1\o h)\ast S'(B((1\o g)_{(3)}))\\
&=&(1\o g_{(1)})\ast B((1\o g_{(2)}))\ast(1\o h)\ast S'(B((1\o g_{(3)})))\\
&=&(1\o g_{(1)})\ast (g_{(2)}\o 1)\ast(1\o h)\ast S'(g_{(3)}\o 1)\\
&=&( g_{(2)}\o (g_{(1)})\ast(1\o h)\ast S'(g_{(3)}\o 1)\\
&=&(g_{(21)}\o g_{(1)}(g_{(22)}\rightharpoonup h))\ast(T(g_{(31)})\o T(g_{(32)})\circ g_{(33)})\\
&=&(g_{(2)}\o g_{(1)}(g_{(3)}\rightharpoonup h))\ast(T(g_{(4)})\o 1)\\
&=&1\o g_{(1)}(g_{(2)}\rightharpoonup h)\\
\hspace{6cm} &=&\psi(g\circ h).
\end{eqnarray*}
\smallskip

\noindent
\eqref{it:3.5b} Define
\begin{center}
$\mathcal{T}=\left\{t_{(H,\c,\circ)}=\psi:I((H,\c,\circ))\rightarrow FJ((H,\c,\circ))\mid \text{for~any~Hopf~brace\ }~(H,\c,\circ)\right\}$,
\end{center}
where $\psi$ is defined by \eqref{it:3.5a}.
We next prove that $\mathcal{T}$ is a natural transformation.

Indeed, for any Hopf brace homomorphism $f: (H,\c,\circ)\rightarrow(H',\c,\circ)$, and for any $h\in (H,\c,\circ)$, we have
$$
FJ(f)\mathcal{T}_{(H,\c,\circ)}(h)=(f\o f)(1\o h)=f(1)\o f(h)=1\o f(h),
$$
and
$$
\mathcal{T}_{(H',\c,\circ)}I(f)(h)=\mathcal{T}_{(H',\c,\circ)}f(h)=1\o f(h),
$$
that is, $FJ(f)\mathcal{T}_{(H,\c,\circ)}=\mathcal{T}_{(H',\c,\circ)}I(f)$. Hence $\mathcal{T}:I\rightarrow FJ$ is a natural transformation.\hspace{7cm}\qedhere
\end{proof}

\begin{exam}
Let $H$ and $K$ be cocommutative Hopf algebras, and $H$ a left $K$-module bialgebra with the module action $\triangleright$. Then the smash product of cocommutative Hopf algebras $G=H\#K$ is a cocommutative Hopf brace (see \cite[Example 1.5]{AGV}). Here
\begin{eqnarray*}
(h\#k)(h'\#k')&=&hh'\#kk', \\
(h\#k)\circ(h'\#k') &=&h(k_{(1)}\triangleright h')\#k_{(2)}k',\\
S_{H\#K}(h\#k)&=&S_H(h)\#S_K(k),\\
T_{H\#K}(h\#k)&=&S_K(k_{(1)})\triangleright S_H(h)\#S_K(k_{(2)}), \ h, h'\in H, k, k'\in K.
\end{eqnarray*}

In the following, applying Theorem~\mref{thm:embed}, we construct a Rota-Baxter Hopf algebra containing $G$. To define the operations on $G':=G\o G=H\#K\o H\#K$: for  $h\#k\o h'\#k',g\#l\o g'\#l'\in G'$, take
\begin{eqnarray*}
(h\#k\o h'\#k')\ast(g\#l\o g'\#l')&:=&h(k_{(1)}\triangleright g)\# k_{(2)}l\o h'(k_{(3)}\triangleright g')\#k'l',\\
S'(h\#k\o h'\#k')&:=&(S_K(k_{(1)})\triangleright S_H(h_{(1)})\#S_K(k_{(2)})\o S_K(k_{(3)})\triangleright S_H(h')\#S_K(k')).
\end{eqnarray*}

Consider the sub-Hopf algebras
$$H:=\{((h\#k)_{(1)}\o (h\#k)_{(2)})\mid h\#k\in G\}, \quad L:=\{((h\#k)\o 1\#1)\mid h\#k\in G\}$$
of $G'$. We obtain the decomposition $G'=H\ast L$, since
\begin{eqnarray*}
h\#k\o h'\#k'&=&((h'\#k')_{(1)}\o (h'\#k')_{(2)})\ast (T((h'\#k')_{(3)})\circ (h\#k)\o 1\#1)\\
&=&(h'_{(1)}\#k'_{(1)}\o h'_{(2)}\#k'_{(2)})\ast ((S_K(k'_{(3)})\triangleright S_H(h'_{(3)}))(S_K(k'_{(4)})\triangleright h)\# S_K(k'_{(5)})k\o 1\#1).
\end{eqnarray*}

Define a Rota-Baxter operator $B:G'\rightarrow G'$ by
\begin{eqnarray*}
B((h\#k,h'\#k'))&=&T(h\#k)\circ (h'\#k')\o 1\#1\\
&=&(S_K(k_{(1)})\triangleright S_H(h))(S_K(k_{(2)})\triangleright h')\# S_K(k_{(3)})k'\o 1\#1.
\end{eqnarray*}
Hence, $G=H\#K$ embeds into the Rota-Baxter Hopf algebra $(G',\ast,B)$ by sending $h\#k$ to $1\#1\o h\#k$.
\mlabel{exam:2.9}
\end{exam}

\begin{remark}
\begin{enumerate}
\item Let $(G, B)$ be a Rota-Baxter group. Then by \cite[Theorem 1]{Go}, we naturally obtain a Rota-Baxter operator $B$ on the group algebra $F[G]$. If $f:G\rightarrow G$ is a group isomorphism, then $\varphi:F[G]\rightarrow F[G]$ is a bialgebra automorphism, where $\varphi(x)=\sum\a_if(g_i)$, for all $x=\sum\a_ig_i\in F(G)$.

 Theorem \mref{prop:rboHb} shows that $(F[G],\cdot,\circ_B)$, $(F[G],\cdot,\circ_{\widetilde{B}})$ and $(F[G],\cdot,\circ_{B^{(\varphi)}})$ are cocommutative Hopf braces. Hence, by Theorem \mref{thm:embed}, they can be embedded into Rota-Baxter Hopf algebras.

 \item Let $L$ be a Lie algebra and $R$ be a Rota-Baxter operator of weight 1 on $L$. Then by \cite[Lemmas 2 and 3]{Go}, the map $R$ can be extended to a linear map $B$: $U(L)\rightarrow U(L)$ such that $B$ is a Rota-Baxter operator on $U(L)$, where $U(L)$ is the universal enveloping algebra of the Lie algebra $L$.

If there is an automorphism $f$ on a Lie algebra $L$, then we naturally have a bialgebra automorphism $\psi:U(L)\rightarrow U(L)$.
 Theorem \mref{prop:rboHb} again shows that $(U(L),\cdot,\circ_B)$, $(U(L),\cdot,\circ_{\widetilde{B}})$ and $(U(L),\cdot,\circ_{B^{(\psi)}})$ are cocommutative Hopf braces. Hence by Theorem \mref{thm:embed}, they can be embedded into Rota-Baxter Hopf algebras.
\end{enumerate}
\mlabel{rem:2.10}
\end{remark}

In what follows, we establish a relationship between cocommutative post-Hopf algebras and Rota-Baxter cocommutative Hopf algebras, which goes beyond the one given in~\mcite{LST}.

 \begin{defn}
 \mcite{LST} A {\bf post-Hopf algebra} is a pair $(H,\blacktriangleright)$, where $H$ is a Hopf algebra and $\blacktriangleright:H\o H\rightarrow H$ is a coalgebra homomorphism satisfying
\begin{eqnarray*}
&&x\blacktriangleright (yz)=(x_{(1)}\blacktriangleright y) (x_{(2)}\blacktriangleright z),\\
&&x\blacktriangleright (y\blacktriangleright z)=(x_{(1)} (x_{(2)}\blacktriangleright y))\blacktriangleright z,\quad x, y, z\in H,
\end{eqnarray*}
and the left multiplication $\a_{\blacktriangleright}: H\rightarrow$ End$(H)$ defined by
$$
\a_{\blacktriangleright,x}(y)=x\blacktriangleright y, ~\text{for all}~ x,y\in H,
$$
is convolutionally invertible in Hom$(H, $ End$(H))$, that is, there exists a unique map $\b_{\blacktriangleright}:H\rightarrow$ End$(H)$ such that
$$
\a_{\blacktriangleright,x_{(1)}}\b_{\blacktriangleright,x_{(2)}}=\b_{\blacktriangleright,x_{(1)}}
\a_{\blacktriangleright,x_{(2)}}=\v(x)\text{id}_H,\  x\in H.
$$
\mlabel{defn:pHA}
\end{defn}

A homomorphism from a post-Hopf algebra $(H, \blacktriangleright)$ to another post-Hopf algebra $(H',\blacktriangleright')$ is a Hopf algebra homomorphism $g : H\rightarrow H'$ satisfying
$$
g(x \blacktriangleright y) = g(x) \blacktriangleright' g(y), ~x, y \in H.
$$

\cite[Theorem 2.13]{LST} states that a cocommutative post-Hopf algebra $(H,\blacktriangleright, S)$ and its subadjacent Hopf algebra $(H,\ast_{\blacktriangleright}, S_{\blacktriangleright})$ form a Hopf brace, where the multiplication $\ast_{\blacktriangleright}$ and the antipode $S_{\blacktriangleright}$ are given by
\begin{equation}
x \ast_{\blacktriangleright} y := x_{(1)} (x_{(2)} \blacktriangleright y), \quad
S_\blacktriangleright(x) := \b_{\blacktriangleright,x_{(1)}}(S (x_{(2)})), \quad x, y\in H.
\mlabel{eq:phhb}
\end{equation}
Conversely, any cocommutative Hopf brace $(H, \c, \circ)$ gives a post-Hopf algebra $(H, \blacktriangleright)$ with the multiplication $\blacktriangleright$ defined
by
\begin{equation}
x \blacktriangleright y=S (x_{(1)}) (x_{(2)} \circ y), x, y \in H.
\mlabel{eq:hbph}
\end{equation}

For a post-Hopf algebra $(H,\blacktriangleright)$ and  $H_\blacktriangleright$ the subadjacent Hopf algebra, by~\cite[Prop.  3.2]{LST}, the identity map $\id_H:H\to H_\blacktriangleright$ is a relative Rota-Baxter operator. We next apply Theorem~\mref{prop:rboHb} to give other avenues to obtain (non-relative) Rota-Baxter operators on Hopf algebras from cocommutative post-Hopf algebras.
Applying Theorems~\mref{prop:rboHb}, \mref{thm:embed} and above relationship between post-Hopf algebras and Hopf braces gives the following correspondences.

\begin{coro}
\begin{enumerate}
\item Let $(H,B)$ be a Rota-Baxter cocommutative Hopf algebra $(H,B)$. Let $(H,\c,\circ_B)$ be the Hopf brace from Theorem~\mref{prop:rboHb}, where $$g\circ_B h:=g_{(1)}B(g_{(2)})hS(B(g_{(3)})),\ g,h\in H.$$
Then, for
$$x\blacktriangleright y := S(x_{(1)})(x_{(2)} \circ_B y)=B(x_{(1)})yS(B(x_{(2)})), \ x,y\in H,$$
the pair $(H, \blacktriangleright)$ is a post-Hopf algebra.
\mlabel{it:3.9a}
\item Let $(H,\blacktriangleright, S)$ be a cocommutative post-Hopf algebra and let $(H,\c,\ast_{\blacktriangleright})$ be the corresponding Hopf brace defined in Eq.~\meqref{eq:phhb}. Define an operator
$$B:H\o H\to H\o H, \  B(x\o y):=S_{\blacktriangleright}(x)\ast_{\blacktriangleright} y\o 1,\ x,y\in H.$$
Then the pair $(H\o H,B)$ is a cocommutative Rota-Baxter Hopf algebra.
\mlabel{it:3.9b}
\item Every cocommutative post-Hopf algebra is embedded as the post-Hopf subalgebra of a Rota-Baxter cocommutative Hopf algebra.
\mlabel{it:3.9c}
     \end{enumerate}
\mlabel{prop:rbopHA}
\end{coro}

\begin{proof} \eqref{it:3.9a} This follows from Theorem \mref{prop:rboHb} and above fact that Eq~\meqref{eq:hbph} defines a post-Hopf algebra.

\eqref{it:3.9b} This follows from the above fact that  $(H,\c,\ast_{\blacktriangleright})$ defined by Eq.~\meqref{eq:phhb} is a Hopf brace together with Theorem \mref{thm:embed} that  $(H\o H, B)$ is a cocommutative Rota-Baxter  Hopf algebra.

\eqref{it:3.9c} Let $H$ be a cocommutative post-Hopf algebra. Then by \eqref{it:3.9b}, $H\o H$ is a cocommutative Rota-Baxter Hopf algebra. Then by \eqref{it:3.9a}, $(H\o H, \blacktriangleright_{H\o H})$ is a post-Hopf algebra with
$$(x\o y) \blacktriangleright_{H\o H} (h\o g) =B(x_{(1)}\o y_{(1)})\ast(h\o g) \ast S'(B(x_{(2)}\o y_{(2)})), \quad x,y,h,g\in H,$$
where $\ast$ and $S'$ are defined in Theorem \mref{thm:embed}.

Embed $H$ into $H\o H$ by the map
$$\psi:H\to H\o H, \ g\mapsto  1\o g.$$
Then $\psi$ is an isomorphism from the post-Hopf algebra $G$ to Im$(\psi)$, regarded as a post-Hopf subalgebra of $(H\o H, \blacktriangleright_{H\o H})$.
This is because, for all $g, h\in H$,  on the one hand,
\begin{eqnarray*}
\psi(g)\blacktriangleright \psi(h)&=& (1\o g)\blacktriangleright (1\o h)\\
&=& B(1\o g_{(1)})\ast(1\o h)\ast S'(B(1\o g_{(2)}))\\
&=& (g_{(1)}\o 1)\ast (1\o h)\ast S'(g_{(2)}\o 1)\\
&=& (g_{(1)}\o S(g_{(2)})(g_{(3)}\ast_{\blacktriangleright} h))\ast (S_{\blacktriangleright}(g_{(2)})\o 1)\\
&=& g_{(1)} \ast_{\blacktriangleright} S_{\blacktriangleright}(g_{(2)})\o S(g_{(3)})(g_{(4)}\ast_{\blacktriangleright} h)\\
&=& 1\o S(h_{(1)})(h_{(2)}\ast_{\blacktriangleright} g)\\
&=& 1\o S(g_{(1)})g_{(2)}(g_{(3)})\blacktriangleright h)\\
&=& 1\o (g\blacktriangleright h)\\
&=&\psi(g\blacktriangleright h),
\end{eqnarray*}
and on the other hand, by Theorem \mref{thm:embed}, $\psi(gh)=\psi(g)\ast\psi(h)=\psi(gh)$.
\end{proof}

\begin{defn}
\mlabel{defn:YBE}
A solution of the {\bf Yang-Baxter equation} on a vector space $V$ is a linear bijection $R:V\o V\rightarrow V\o V$ such that
$$
(R\o \text{id})(\text{id}\o R)(R\o \text{id})=(\text{id}\o R)(R\o \text{id})(\text{id}\o R).
$$
\end{defn}

In the following, we obtain solutions of the Yang-Baxter equation via Rota-Baxter operators on cocommutative Hopf algebras.

Let $(H,B)$ be a Rota-Baxter cocommutative Hopf algebra. By Theorem~\ref{prop:rboHb}\eqref{it:3.1a},  $(H,\c,\circ_{B})$ is a cocommutative Hopf brace. By \cite[Lemma 1.8 and Lemma 2.2]{AGV}, $H$ is a left $H_B$ module bialgebra via $\rightharpoonup_B$, and $H$ is a right $H_B$ module coalgebra via $\leftharpoonup_B$, where the actions $\rightharpoonup_B$ and $\leftharpoonup_B$ are defined by
	\begin{equation}
		h\rightharpoonup_B k:=S(h_{(1)})(h_{(2)}\circ_B k), ~~~h\leftharpoonup_B k:=T(h_{(1)}\rightharpoonup k_{(1)})\circ_B h_{(2)}\circ_Bk_{(2)}, \ h, k\in H,
		\mlabel{eq:arrowb}
	\end{equation}
for $T$ in Eq.~\meqref{eq:desc}.

\begin{prop}
Let $(H,B)$ be a Rota-Baxter cocommutative Hopf algebra. Then the linear map
\begin{equation}
c:H\o H\rightarrow H\o H, \quad c(x\o y)=x_{(1)}\rightharpoonup_B y_{(1)}\o x_{(2)}\leftharpoonup_B y_{(2)},
\label{eq:ybebrace}
\end{equation}
is a coalgebra isomorphism and a solution of the Yang-Baxter equation on the vector space $H$.
Furthermore, the actions $\rightharpoonup_B$ and $\leftharpoonup_B$ are given in terms of $B$ by
\begin{eqnarray}
h\rightharpoonup_B k&=&B(h_{(1)})kS(B(h_{(2)})),\  \label{eq:arrowr}\\
h\leftharpoonup_B k&=&S(B(h_{(1)}\rightharpoonup k_{(1)}))S(h_{(2)}\rightharpoonup k_{(2)})h_{(3)}(h_{(4)}\rightharpoonup k_{(3)})B(h_{(5)}\rightharpoonup k_{(4)}),\ \  h, k\in H,
\label{eq:arrowl}
\end{eqnarray}
respectively.
\mlabel{prop:rboYBE}
\end{prop}

\begin{proof} Since  $(H,B)$ is a Rota-Baxter cocommutative Hopf algebra, by Theorem~\ref{prop:rboHb}\eqref{it:3.1a}, $(H,\c,\circ_{B})$ is a Hopf brace. By \cite[Corollary 2.4]{AGV}, the Hopf brace $(H,\c,\circ_{B})$ gives a solution of the Yang-Baxter equation by Eq.~\meqref{eq:ybebrace}.
	
Furthermore, Eqs.~\eqref{eq:arrowr} and \eqref{eq:arrowl} are verified as follows.
\begin{eqnarray*}
	h\rightharpoonup_B k&=&S(h_{(1)})(h_{(2)}\circ_B k)\\
	&=&S(h_{(1)})h_{(2)(1)}B(h_{(2)(2)})kS(B(h_{(2)(3)}))\\
	&=&B(h_{(1)})kS(B(h_{(2)})).
\end{eqnarray*}
By \cite[Prop.  5 and Corollary 3]{Go} and Lemma~\mref{lemma:desHA}, we obtain
\begin{eqnarray*}
	h\leftharpoonup_B k&=&T(h_{(1)}\rightharpoonup k_{(1)})\circ_B h_{(2)}\circ_Bk_{(2)}\\
	&=&T(h_{(1)}\rightharpoonup k_{(1)})\circ_B(h_{(2)}B(h_{(3)})k_{(2)}S(B(h_{(4)})))\\
	&=&T(h_{(1)}\rightharpoonup k_{(1)})B(T(h_{(2)}\rightharpoonup k_{(2)}))h_{(3)}B(h_{(4)})k_{(3)}S(B(h_{(5)}))S(B(T(h_{(6)}\rightharpoonup k_{(4)})))\\
	&=&S(B(h_{(1)}\rightharpoonup k_{(1)}))S(h_{(2)}\rightharpoonup k_{(2)})B(h_{(3)}\rightharpoonup k_{(3)})\\
	&&\times B(S(B(h_{(4)}\rightharpoonup k_{(4)}))S(h_{(5)}\rightharpoonup k_{(5)})B(h_{(6)}\rightharpoonup k_{(6)}))\\
	&&\times h_{(7)}(h_{(8)}\rightharpoonup k_{(7)})S(B(T(h_{(9)}\rightharpoonup k_{(8)})))\\
	&=&S(B(h_{(1)}\rightharpoonup k_{(1)}))S(h_{(2)}\rightharpoonup k_{(2)})\underline{B(h_{(3)}\rightharpoonup k_{(3)})B(T(h_{(4)}\rightharpoonup k_{(4)}))}\\
	&&\times h_{(5)}(h_{(6)}\rightharpoonup k_{(5)})S(B(T(h_{(7)}\rightharpoonup k_{(6)})))\\
	&=&S(B(h_{(1)}\rightharpoonup k_{(1)}))S(h_{(2)}\rightharpoonup k_{(2)})h_{(3)}(h_{(4)}\rightharpoonup k_{(3)})S(B(T((h_{(5)}\rightharpoonup k_{(4)}))))\\
	&& \qquad (\text{applying Lemma \mref{lemma:desHA}\eqref{lemma:rboproperty3} to the underlined part})\\
	&=&S(B(h_{(1)}\rightharpoonup k_{(1)}))S(h_{(2)}\rightharpoonup k_{(2)})h_{(3)}(h_{(4)}\rightharpoonup k_{(3)})B(T(T((h_{(5)}\rightharpoonup k_{(4)}))))\\
	&=&S(B(h_{(1)}\rightharpoonup k_{(1)}))S(h_{(2)}\rightharpoonup k_{(2)})h_{(3)}(h_{(4)}\rightharpoonup k_{(3)})B(h_{(5)}\rightharpoonup k_{(4)}).  \hspace{3cm} \qedhere
\end{eqnarray*}
\end{proof}

\subsection{Relative Rota-Baxter operators on cocommutative Hopf algebras, Hopf braces and bijective 1-cocycles}\
In this subsection, we construct a Hopf brace from a relative Rota-Baxter operator on Hopf algebra and establish the connection between relative Rota-Baxter operators and bijective 1-cocycles on Hopf algebras.

In \mcite{LST}, the authors generalized the notion of Rota-Baxter operators on a Hopf algebra~\mcite{Go} to that of relative Rota-Baxter operators, by generalizing the adjoint action to arbitrary actions on bialgebras.

\begin{defn}
\mlabel{defn:rRbo}
\mcite{LST} Let $H$ and $K$ be Hopf algebras such that $K$ is a left $H$-module bialgebra via an action $\rightharpoonup$. A coalgebra homomorphism $\tau:K\rightarrow H$ is called a {\bf relative Rota-Baxter operator} with respect to the left $H$-module bialgebra $(K,\rightharpoonup)$ if
\begin{equation}
\tau(a)\tau(b)=\tau(a_{(1)}(\tau(a_{(2)})\rightharpoonup b)), \quad a, b\in K.
\mlabel{eq:relrbo}
\end{equation}
\end{defn}

\begin{lemma}\mcite{LST}\mlabel{lem:rRbodes}
Let $\tau:K\rightarrow H$ be a relative Rota-Baxter operator with respect to a cocommutative $H$-module bialgebra $(K,\rightharpoonup)$. Then $(K,\circ,\D,S_{\tau})$ is a Hopf algebra, which is called the {\bf descendant Hopf algebra} and denoted by $K_{\tau}$, where the antipode $S_{\tau}$ and the multiplication $\circ$ is given by
$$S_{\tau}(k)=S_H(\tau(k_{(1)}))\rightharpoonup S_K(k_{(2)}),$$
and
$$
k\circ k'=k_{(1)}(\tau (k_{(2)})\rightharpoonup k'),~ k,k'\in K.
$$
\end{lemma}

In the following, we will prove that the descendant Hopf algebra also forms a Hopf brace. In fact, The descendant Hopf algebra derived from a relative Rota-Baxter operator generalizes the descendant Hopf algebra obtained through a Rota-Baxter Hopf algebra in Lemma \mref{lemma:desHA}.

\begin{prop}
Let $\tau:K\rightarrow H$ be a relative Rota-Baxter operator with respect to a cocommutative $H$-module bialgebra $(K,\rightharpoonup)$. Then $(K,K_{\tau},\D)$ is a Hopf brace.

Conversely, let $(H, \c, \circ)$ be a cocommutative Hopf brace. Then with the action $\rightharpoonup$ defined by Lemma \mref{lemma:Hbproperty}, the triple $(H, H_{\circ}, \rightharpoonup,\id_H)$ is a relative Rota-Baxter Hopf algebra.
\mlabel{prop:rRboHb}
\mlabel{prop:HbrRbo}
\end{prop}

\begin{proof}
Indeed, for any $k\in K$, we have
\begin{eqnarray*}
(k_{(1)}\circ k')S(k_{(2)})(k_{(3)}\circ k'')&=&k_{(1)}(\tau (k_{(2)})\rightharpoonup k')S(k_{(3)})k_{(4)}(\tau (k_{(5)})\rightharpoonup k'')\\
&=&k_{(1)}(\tau (k_{(2)})\rightharpoonup k')(\tau (k_{(3)})\rightharpoonup k'')\\
&=&k_{(1)}(\tau (k_{(2)})\rightharpoonup (k'k''))\\
 &=& k\circ (k'k'').
\end{eqnarray*}

For the second claim, taking $\tau=\id_H$ in Eq.~\meqref{eq:relrbo} yields Eq. \meqref{eq:hbm}.
\end{proof}

We next recall the notion of bijective 1-cocycles on Hopf algebras as a special case of 1-cocycles, also known as crossed homomorphisms~\mcite{Sw} and recently identified as relative difference operators~\mcite{GLST}.

\begin{defn}
\mlabel{defn:cocycle}
\mcite{AGV} Let $H$ and $A$ be Hopf algebras, and $A$ a left $H$-module algebra via an action $\rightharpoonup$. A {\bf bijective 1-cocycle} is a coalgebra isomorphism $\pi : H\rightarrow A$ such that
$$
\pi(hk)=\pi(h_{(1)})(h_{(2)} \rightharpoonup \pi(k)), \quad h, k\in H.
$$
\end{defn}

Bijective 1-cocycles can be related to relative Rota-Baxter operators in multiple ways. We first give a direct connection.

\begin{prop} Let $H$ and $A$ be Hopf algebras, and $A$ an $H$-module bialgebra via an action $\rightharpoonup$. If $\tau:A\rightarrow H$ is a invertible linear map, then $\tau$ is a relative Rota-Baxter operator if and only if $\tau^{-1}$ is a bijective 1-cocycle.
\mlabel{prop:rRbococycle}
\end{prop}

\begin{proof} If $\tau:A\rightarrow H$ is a relative Rota-Baxter operator with respect to the left $H$-module bialgebra $(A,\rightharpoonup)$, we have
$$\tau(a)\tau(b)=\tau(a_{(1)}(\tau(a_{(2)})\rightharpoonup b )),
$$
for all $a,b\in A$. By the bijectivity of $\tau$, there exist $h,k\in H$ such that $\tau^{-1}(h)=a$ and $\tau^{-1}(k)=b$. Then $hk=\tau(\tau^{-1}(h_{(1)})(h_{(2)}\rightharpoonup \tau^{-1}(k)))$, that is, $$\tau^{-1}(hk)=\tau^{-1}(h_{(1)})(h_{(2)}\rightharpoonup \tau^{-1}(k)).
$$
Moreover, since $\tau$ is a coalgebra isomorphism, it is easy to verify that $\tau^{-1}$ is a coalgebra isomorphism. Therefore, $\tau^{-1}$ is a bijective 1-cocycle.

Conversely, if $\tau^{-1}:H\rightarrow A$ is a bijective 1-cocycle, then we have
$$
\tau^{-1}(hk)=\tau^{-1}(h_{(1)})(h_{(2)} \rightharpoonup \tau^{-1}(k)),\ h, h\in H,
$$
and there exist $a,b\in A$ such that $\tau(a)=h$ and $\tau(b)=k$.
Hence $\tau^{-1}(\tau(a)\tau(b))=a_{(1)}(\tau(a_{(2)}) \rightharpoonup b)$ which means that
$$\tau(a)\tau(b)=\tau(a_{(1)}(\tau(a_{(2)})\rightharpoonup b )).
$$
Hence $\tau$ is a relative Rota-Baxter operator.
\end{proof}

Applying the relationship between Hopf braces and bijective 1-cocycles, we arrive at the next connection between 1-cocycles and relative Rota-Baxter operators on Hopf algebras.

\begin{prop} Let $H$ and $K$ be cocommutative Hopf algebras.
\begin{enumerate}	
\item If $\pi:H\rightarrow K$ is a bijective 1-cocycle, then id$_K:K\rightarrow K_{\circ_{\pi}}$ is a relative Rota-Baxter operator, where $\circ_{\pi}$ is defined by $a\circ_{\pi} b=\pi(\pi^{-1}(a)\pi^{-1}(b))$, for all $a, b \in K$.
\mlabel{it:3.21a}
\item If $\tau:K\rightarrow H$ is a relative Rota-Baxter operator with respect to a left $H$-module bialgebra $(K,\rightharpoonup)$, then  id$_K:K\rightarrow K_{\circ}$ is a relative Rota-Baxter operator, where $\circ$ is defined by $x\circ y=x_{(1)}(\tau(x_{(2)})\rightharpoonup y)$, for all $x,y\in K$.
\mlabel{it:3.21b}
\end{enumerate}
\mlabel{prop:cocyclerRbo}
\end{prop}

\begin{proof}
\eqref{it:3.21a} By \cite[Theorem 1.12]{AGV}, if $\pi:H\rightarrow K$ is a bijective 1-cocycle, then $(K,\c, \circ_{\pi})$ is a Hopf brace, with $a\circ_{\pi} b=\pi(\pi^{-1}(a)\pi^{-1}(b))$ for $a, b \in K$. So by Lemma \mref{lemma:Hbproperty},  $K$ is a left $K_{\circ_{\pi}}$-module bialgebra and $a\circ_\pi b=a_{(1)}(a_{(2)}\rightharpoonup b)$ for all $a,b\in K$ means that  id$_K:K\rightarrow K_{\circ_{\pi}}$ is a relative Rota-Baxter operator.
\smallskip

\noindent
\eqref{it:3.21b} If $\tau:K\rightarrow H$ is a relative Rota-Baxter operator with respect to a left $H$-module bialgebra $(K,\rightharpoonup)$, then by Proposition \mref{prop:rRboHb}, there exists a Hopf brace $(K,\c,\circ)$ with $x\circ y=x_{(1)}(\tau(x_{(2)})\rightharpoonup y)$ for $x,y\in K$. So by Definition \mref{defn:cocycle} and Lemma \mref{lemma:Hbproperty}, id$_K: K_{\circ}\rightarrow K$ is a bijective 1-cocycle and its inverse is a relative Rota-Baxter operator.
\end{proof}

\begin{coro}
Let $(H,B)$ be a Rota-Baxter cocommutative Hopf algebra. Then  id$_H: H\rightarrow H_B$ is a relative Rota-Baxter operator.
\mlabel{coro:rborRbo}
\end{coro}
\begin{proof} By Theorem \mref{prop:rboHb}, $(H,\c,\circ_B)$ is a Hopf brace. So $H$ is a left $H_B$-module bialgebra with the action $\rightharpoonup$ such that $a \circ_B b = a_{(1)}(a_{(2)} \rightharpoonup  b)$, for $a,b\in H$. Hence, by the above two definitions, the identity map id$_H: H\rightarrow H_B$ is a relative Rota-Baxter operator and id$_H:H_B\rightarrow H$ a bijective 1-cocycle.
\end{proof}

\begin{remark}
\mlabel{remark:cocyclerbo}
Let $\pi : H\rightarrow A$ be a bijective 1-cocycle. By \cite[Theorem 1.12]{AGV} and Theorem \mref{thm:embed}, we get a Rota-Baxter cocommutative Hopf algebra $(A\o A,B)$, with the Rota-Baxter operator $B(x\o y):=\pi(S\pi^{-1}(x)\pi^{-1}(y))\o 1$, for all $x,y\in A$, and tensor coalgebra $A\o A$, and the multiplication $\ast$ and antipode $S'$ defined by
$$
(x\o y)\ast(z\o t)=\pi(\pi^{-1}(x_{(1)})\pi^{-1}(z))\o y(S(x_{(2)})(\pi(\pi^{-1}(x_{(3)})\pi^{-1}(t)))),
$$
$$
S'(x\o y)=\pi S\pi^{-1}(x_{(1)})\o \pi(S\pi^{-1}(x_{(2)})\pi^{-1}(x_{(3)}S(y))).
$$
\end{remark}

\section{Rota-Baxter operators and symmetric Hopf braces}
\mlabel{sec:shopf}
In this section, we introduce the notion of a symmetric Hopf brace which generalizes the symmetric skew brace (See\mcite{BNY}), and give their characterizations through Rota-Baxter operators of Hopf algebras.

\begin{defn}
A Hopf brace $(H,\cdot,\circ)$ is called {\bf symmetric} if $(H,\circ,\cdot)$ is also a Hopf brace.
\mlabel{defn:symmetric}
\end{defn}

We first characterize symmetric cocommutative Hopf braces by equations.

\begin{prop}
\mlabel{prop:sHb}
A cocommutative Hopf brace $(H,\cdot,\circ)$ is symmetric if and only if  the action $\rightharpoonup$ in Lemma \mref{lemma:Hbproperty} satisfies
$$
ab_{(1)}(b_{(2)}\rightharpoonup c)=a_{(1)}b_{(1)}((a_{(2)}b_{(2)})\rightharpoonup (T(a_{(3)})\rightharpoonup c)),
\ a,b,c\in H.$$
\end{prop}

\begin{proof}
If a cocommutative Hopf brace $(H,\cdot,\circ)$ is symmetric, then Lemma \mref{lemma:Hbproperty} holds and $(H,\circ,\c)$ is a Hopf brace, that is, $a(b\circ c)=(a_{(1)}b)\circ T(a_{(2)})\circ (a_{(3)}c)$, for all $a,b,c\in H$.
Hence
\begin{eqnarray*}
ab_{(1)}(b_{(2)}\rightharpoonup c)&\stackrel{\eqref{eq:hbm}}{=}&a(b\circ c)\\
&=&(a_{(1)}b)\circ T(a_{(2)})\circ (a_{(3)}c)\\
&=&(a_{(1)}b_{(1)})((a_{(2)}b_{(2)})\rightharpoonup(\underline{T(a_{(3)})(T(a_{(4)})
	\rightharpoonup a_{(5)})}(T(a_{(6)})\rightharpoonup c)))\\	&\stackrel{\eqref{eq:hbm}}{=}&(a_{(1)}b_{(1)})((a_{(2)}b_{(2)})\rightharpoonup((T(a_{(3)})\circ a_{(4)})(T(a_{(5)})\rightharpoonup c)))\\
	&=&a_{(1)}b_{(1)}((a_{(2)}b_{(2)})\rightharpoonup (T(a_{(3)})\rightharpoonup c)).
\end{eqnarray*}
Here the underlined part indicates where a simplification occurs.

Conversely, if the action $\rightharpoonup$ satisfies the given equation, we prove that $(H,\circ,\c)$ is also a Hopf brace, that is, for all $a,b,c\in H$,
\begin{eqnarray*}
&&a(b\circ c)=(a_{(1)}b)\circ T(a_{(2)})\circ (a_{(3)}c).
\end{eqnarray*}
Because $(H,\cdot,\circ)$ is a cocommutative Hopf brace, Lemma \mref{lemma:Hbproperty} holds.

Hence, for any $a,b,c\in H$,
\begin{eqnarray*}
a(b\circ c)&=&a b_{(1)}(b_{(2)}\rightharpoonup c)\\
&=&a_{(1)}b_{(1)}((a_{(2)}b_{(2)})\rightharpoonup (T(a_{(3)})\rightharpoonup c))\\
&=&(a_{(1)}b)\circ T(a_{(2)})\circ (a_{(3)}c),
\end{eqnarray*}
where the last step follows from the proof of necessity above.
\end{proof}

We next study symmetric Hopf braces by means of left $H^{\rm op}$-modules.
\begin{defn}
\mlabel{defn:moduleHb}
Let $(H,\cdot,\circ)$ be a Hopf brace. Let $H^{\rm op}$ denote the opposite Hopf algebra of $(H,\cdot)$.
We call $(H,\cdot,\circ)$ an {\bf $H^{\rm op}$-module Hopf brace} if $(H, \rightharpoonup)$ is a left $H^{\rm op}$-module, where the action $\rightharpoonup$ is defined by
\begin{equation}
a\rightharpoonup b:=S(a_{(1)})(a_{(2)}\circ b), \quad a, b\in H.
\mlabel{eq:hbprod}
\end{equation}
\end{defn}

Note that the $H^{\rm op}$-module $(H, \rightharpoonup)$ in Definition \mref{defn:moduleHb} is actually a left $H^{\rm op}$-module bialgebra by Lemma \mref{lemma:Hbproperty}.

Let $(H, \cdot, \circ)$ be an $H^{\rm op}$-module Hopf brace. By the proof of \cite[Lemma 1.8]{AGV}, for all $a, b, c\in H$, we have $(a\circ b)\rightharpoonup c=a\rightharpoonup(b\rightharpoonup c)$. That is to say, $(H, \rightharpoonup)$ is an $H_{\circ}$-module. Since $(H, \rightharpoonup)$ is an $H^{\rm op}$-module, we have $a\rightharpoonup(b\rightharpoonup c)=(b a)\rightharpoonup c$. Therefore,
$$
(a\circ b)\rightharpoonup c=(b a)\rightharpoonup c.
$$

\begin{exam}
	Let $(H,\cdot)$ be a Hopf algebra. Define an operator $\circ$ on the vector space $H$ by
	$$
	a\circ b:=ba, \quad a,b\in H.
	$$
Then $(H,\cdot,\circ)$ is a Hopf brace.
Let the action $\rightharpoonup$ on $H$ be given in Eq.~\meqref{eq:hbprod}. Then for any $a,b\in H,$ we have  $a\rightharpoonup b=S(a_{(1)})(a_{(2)}\circ b)=S(a_{(1)})ba_{(2)}$, so that $(H,\rightharpoonup)$ is a module algebra. Hence every Hopf algebra $(H,\cdot)$ is naturally an $H^{\rm op}$-module Hopf brace $(H,\cdot,\circ)$.
	\mlabel{exam:exHb}
\end{exam}

\begin{prop}
Every cocommutative  $H^{\rm op}$-module Hopf brace $(H,\cdot,\circ)$ is symmetric.
\mlabel{prop:smHb}
\end{prop}
\begin{proof} If $H$ is an $H^{\rm op}$-module Hopf brace, then $(ab)\rightharpoonup c=(b\circ a)\rightharpoonup c$, for all $a,b,c\in H$.

Hence,
\begin{eqnarray*}
a_{(1)}b_{(1)}((a_{(2)}b_{(2)})\rightharpoonup (T(a_{(3)})\rightharpoonup c))&=&a_{(1)}b_{(1)}((b_{(2)}\circ a_{(2)})\rightharpoonup (T(a_{(3)})\rightharpoonup c))\\
&=&a_{(1)}b_{(1)}(b_{(2)}\rightharpoonup (a_{(2)}\rightharpoonup (T(a_{(3)})\rightharpoonup c))\\
&=&a_{(1)}b_{(1)}(b_{(2)}\rightharpoonup ((a_{(2)}\circ T(a_{(3)}))\rightharpoonup c))\\
&=&ab_{(1)}(b_{(2)}\rightharpoonup c).
\end{eqnarray*}
By Proposition \mref{prop:sHb}, $H$ is symmetric.
\end{proof}

We give a natural construction of symmetric Hopf braces.

\begin{defn}
A cocommutative Hopf algebra $H$ is said {\bf to allow a factorization}\footnote{This terminology is used to avoid confusion with the factorizable Hopf algebra defined by Reshetikhin and Semenov-Tian-Shansky~\mcite{RS}.} if there are nontrivial Hopf subalgebras $A$ and $B$ of $H$ such that $H\cong A\o B$ as a tensor product algebra and $H=AB$ as a Hopf algebra.
Then we denote $H=A\o B$.
\mlabel{de:facthopf}
\end{defn}

Let $H=A\o B$ be a cocommutative Hopf algebra that allows a factorization.
Define a operator $\circ$ on $H$ by
$$
(ab) \circ (a'b'):=aa'b'b, \quad a, a'\in A, b, b'\in B.
$$
Then, it is not difficult to check that $(H, \cdot,\circ)$ is a Hopf brace. It is easy to see that $(H,\circ)$ is isomorphic to the tensor algebra $A\o B^{\rm op}$, where $B^{\rm op}$ is the opposite algebra of $B$.

Further define an action $\rightharpoonup: (H,\circ)\o H\rightarrow H$ by
$$
x\rightharpoonup y:=S(x_{(1)})(x_{(2)}\circ y).
$$
 If $x=ab$, $y=a'b'$, where $a,a'\in A$, $b,b'\in B$, then
 $$
 x\rightharpoonup y=S((ab)_{(1)})((ab)_{(2)}\circ (a'b'))=S(a_{(1)}b_{(1)})(a_{(2)}a'b'b_{(2)})=\v(a)S(b_{(1)})yb_{(2)}.
 $$
Hence,
 $$
 (xy)\rightharpoonup z=\v(aa')S((bb')_{(1)})z(bb')_{(2)}=\v(aa')S(b'_{(1)})S(b_{(1)})zb_{(2)}b'_{(2)}
 =y\rightharpoonup (x\rightharpoonup z).
 $$
Thus we have the following conclusion.
\begin{prop}
For a Hopf algebra $H$ that allows a factorization, the Hopf brace $(H,\cdot,\circ)$ constructed above is an $H^{\rm op}$-module Hopf brace, hence is symmetric.
\mlabel{prop:mHbs}
\end{prop}

We give another construction of symmetric Hopf braces.
\begin{theorem}
Let $(H,\c,\D,S)$ be a cocommutative Hopf algebra, and $(H,\rightharpoonup)$ an $H^{\rm op}$-module bialgebra such that
$$a_{(1)}(a_{(2)}\rightharpoonup b)\rightharpoonup c=(ba)\rightharpoonup c, \quad a, b, c\in H.
$$
Define
$$a\circ b:=a_{(1)}(a_{(2)}\rightharpoonup b) \hbox{ and } T(a):=S(a_{(1)})\rightharpoonup S(a_{(2)})$$
Then $(H,\c, \circ)$ is a symmetric Hopf brace.
\mlabel{thm:coHA}
\end{theorem}

\begin{proof}
We first check that $(H,\circ)$ is an algebra.

This is because, for any $a,b,c\in G$, we have
\begin{eqnarray*}
	(a\circ b)\circ c&=&(a_{(1)}(a_{(2)}\rightharpoonup b))\circ c\\
	&=&(a_{(1)}(a_{(2)}\rightharpoonup b))_{(1)}((a_{(1)}(a_{(2)}\rightharpoonup b))_{(2)}\rightharpoonup c)\\
	&=&a_{(1)}(a_{(2)}\rightharpoonup b_{(1)})((a_{(3)}(a_{(4)}\rightharpoonup b_{(2)}))\rightharpoonup c)\\
	&=&a_{(1)}(a_{(2)}\rightharpoonup b_{(1)})(b_{(2)}a_{(3)}\rightharpoonup c)\\
	&=&a_{(1)}(a_{(2)}\rightharpoonup b_{(1)})(a_{(3)}\rightharpoonup(b_{(2)}\rightharpoonup c))\\
	&=&a_{(1)}(a_{(2)}\rightharpoonup (b_{(1)}(b_{(2)}\rightharpoonup c)))\\
	&=&a_{(1)}(a_{(2)}\rightharpoonup (b\circ c))\\
	&=&a\circ (b\circ c).
\end{eqnarray*}
Hence $(H,\circ)$ is an algebra with unit $1$.

For any $a,b\in H$, we have
\begin{eqnarray*}
\D(a\circ b)&=&\D(a_{(1)}(a_{(2)}\rightharpoonup b))\\
&=&a_{(1)(1)}(a_{(2)}\rightharpoonup b)_{(1)}\o a_{(1)(2)}(a_{(2)}\rightharpoonup b)_{(2)}\\
&=&a_{(1)}(a_{(2)}\rightharpoonup b_{(1)})\o a_{(3)}(a_{(4)}\rightharpoonup b_{(2)})\\
&=&a_{(1)}\circ b_{(1)}\o a_{(2)}\circ b_{(2)}.
\end{eqnarray*}
and $\v(a\circ b)=\v(a_{(1)}(a_{(2)}\rightharpoonup b))=\v(a_{(1)})\v(a_{(2)}\rightharpoonup b)=\v(a_{(1)})\v(a_{(2)})\v(b)=\v(a)\v(b)$.

Hence $(H,\circ,\D)$ is a bialgebra.

In what follows, we prove that $a\circ(bc)=(a_{(1)}\circ b)S(a_{(2)})(a_{(3)}\circ c)$, for any $a,b,c\in H$.
\begin{eqnarray*}
(a_{(1)}\circ b)S(a_{(2)})(a_{(3)}\circ c)
&=&a_{(1)}(a_{(2)}\rightharpoonup b)S(a_{(3)})a_{(4)}(a_{(5)}\rightharpoonup c)\\
&=&a_{(1)}(a_{(2)}\rightharpoonup b)(a_{(3)}\rightharpoonup c)\\
&=&a_{(1)}(a_{(2)}\rightharpoonup (bc))\\
&=&a\circ(bc).
\end{eqnarray*}
 To show that $(H,\circ,\D)$ is a Hopf algebra with the antipode $T(a)=S(a_{(1)})\rightharpoonup S(a_{(2)})$ for $a\in H$. We first prove that $$S(a\rightharpoonup b)=a\rightharpoonup S(b).$$
This is because
\begin{eqnarray*}
S(a\rightharpoonup b)&=&S(S(a_{(1)})(a_{(2)}\circ b))\\
&=&S(a_{(2)}\circ b)a_{(1)}\\
&=&S(a_{(1)}\circ b)a_{(2)}\\
&=&S(a_{(1)}\circ b_{(1)})\v(b_{(2)})a_{(2)}\\
&=&S(a_{(1)}\circ b_{(1)})(\underline{a_{(2)}\circ (b_{(2)}S(b_{(3)}))})\\
&\stackrel{\eqref{eq:brace}}{=}&\underline{S(a_{(1)}\circ b_{(1)})(a_{(2)}\circ b_{(2)})}S(a_{(3)})(a_{(4)}\circ S(b_{(3)}))\\
&=&\v(a_{(1)}\circ b_{(1)})S(a_{(2)})(a_{(3)}\circ S(b_{(2)})) \quad \text{(by antipode property)}\\
&=&S(a_{(1)})(a_{(2)}\circ S(b))\\
&=&a\rightharpoonup S(b).
\end{eqnarray*}

Secondly, we prove $T^2=\id$: for any $a\in H$, we have
\begin{eqnarray*}
T^2(a)&=&T(S(a_{(1)})\rightharpoonup S(a_{(2)}))\\
&=&S(S(a_{(1)})\rightharpoonup S(a_{(2)}))\rightharpoonup S(S(a_{(3)})\rightharpoonup S(a_{(4)}))\\
&=&(S(a_{(1)})\rightharpoonup a_{(2)})\rightharpoonup (S(a_{(3)})\rightharpoonup a_{(4)})\\
&=&(S(a_{(3)})(S(a_{(1)})\rightharpoonup a_{(2)}))\rightharpoonup a_{(4)}\\
&=&(S(a_{(2)})(S(a_{(1)})\rightharpoonup a_{(3)}))\rightharpoonup a_{(4)}\\
&=&(a_{(2)}S(a_{(1)}))\rightharpoonup a_{(3)}\\
&=&\v(a_1)a_2\\
&=&a.
\end{eqnarray*}

In the following, we prove that $T$ is antipode of $(H, \circ,\D)$.

Because $H$ is cocommutative, for any $a\in H$, we get
$$\D(T(a))=T(a_{(1)})\o T(a_{(2)}),$$
and
$$a_{(1)}\rightharpoonup T(a_{(2)})=a_{(1)}\rightharpoonup (S(a_{(2)})\rightharpoonup S(a_{(3)}))=(S(a_{(2)})a_{(1)})\rightharpoonup S(a_{(3)})=S(a).$$
Furthermore, we obtain
$$S(T(a))=T(a)_{(1)}\rightharpoonup T(T(a)_{(2)})=T(a_{(1)}) \rightharpoonup a_{(2)}.$$
Hence
$$
T(a_{(1)})\circ a_{(2)}=T(a_{(1)})(T(a_{(2)})\rightharpoonup a_{(3)})=T(a_{(1)})S(T(a_{(2)}))=\v(a)1,
$$
$$
a_{(1)}\circ T(a_{(2)})=a_{(1)}(a_{(2)}\rightharpoonup T(a_{(3)}))=a_{(1)}S(a_{(2)})=\v(a)1.
$$

 Now we have shown that $(H,\cdot,\circ)$ is a Hopf brace, and an $H^{\rm op}$-module Hopf brace, and therefore a symmetric Hopf brace by Proposition \mref{prop:smHb}.
\end{proof}

\begin{prop}
\mlabel{prop:sHbrbo}
\mlabel{prop:rbosHb}
Let $(H,B)$ be a Rota-Baxter cocommutative Hopf algebra. Let $\vartriangleright$ be the  adjoint action of $H$ on itself given by
$a\vartriangleright b:=a_{(1)}bS(a_{2})$.
\begin{enumerate}
\item
\mlabel{it:rbshb1}
The triple $(H,\c,\circ_B)$ is a symmetric Hopf brace if and only if
$$
ab_{(1)}(B(b_{(2)})\vartriangleright c)=a_{(1)}b_{(1)}((B(a_{(2)}b_{(2)})B(T(a_{(3)})))\vartriangleright c), \quad a, b, c\in H.
$$
\item
\mlabel{it:rbshb2}
$(H,\c,\circ_B)$ is an $H^{\rm op}$-module Hopf brace if and only if
$$B(ba)\vartriangleright c=(B(a)B(b))\vartriangleright c, \quad a,b,c\in H.
$$
\item
\mlabel{it:rbshb3}
If $B$ is an anti-algebra homomorphism,
then $(H,\c,\circ_B)$ is a symmetric Hopf brace.
\end{enumerate}
\end{prop}

\begin{proof}
\meqref{it:rbshb1} It is obvious that $(H,\c,\circ_B)$ is a Hopf brace by Theorem \mref{prop:rboHb}. So for the action $\rightharpoonup$ in Definition \mref{defn:moduleHb}: $a\rightharpoonup b=S(a_{(1)})(a_{(2)}\circ_Bb)$ for $a,b\in H$, we have $a\rightharpoonup b=B(a)\rhd b$.

By Proposition \mref{prop:sHb}, A cocommutative Hopf brace $(H,\cdot,\circ)$ is symmetric if and only if
$$
ab_{(1)}(b_{(2)}\rightharpoonup c)=a_{(1)}b_{(1)}((a_{(2)}b_{(2)})\rightharpoonup (T(a_{(3)})\rightharpoonup c)),
\ a,b,c\in H.$$

Note that
\begin{eqnarray*}
a_{(1)}b_{(1)}((a_{(2)}b_{(2)})\rightharpoonup (T(a_{(3)})\rightharpoonup c))&=&a_{(1)}b_{(1)}(B(a_{(2)}b_{(2)})\vartriangleright (B(T(a_{(3)}))\vartriangleright c))\\
&=&a_{(1)}b_{(1)}((B(a_{(2)}b_{(2)})B(T(a_{(3)})))\vartriangleright c),
\end{eqnarray*}
and
\begin{eqnarray*}
ab_{(1)}(b_{(2)}\rightharpoonup c)&=&ab_{(1)}(B(b_{(2)})\vartriangleright c).
\end{eqnarray*}
Hence $(H,\c,\circ_B)$ is a symmetric Hopf brace if and only if
$$
ab_{(1)}(B(b_{(2)})\vartriangleright c)=a_{(1)}b_{(1)}((B(a_{(2)}b_{(2)})B(T(a_{(3)})))\vartriangleright c).
$$

\smallskip

\noindent
\meqref{it:rbshb2}
It is easy to see that $(H,\c,\circ_B)$ is an $H^{\rm op}$-module Hopf brace if and only if $(ba)\rightharpoonup c=a\rightharpoonup(b\rightharpoonup c)$, for all $a,b,c\in H$. Since
$$(ba)\rightharpoonup c=B(ba)\vartriangleright c, \quad a\rightharpoonup(b\rightharpoonup c)=B(a)\vartriangleright (B(b)\vartriangleright c)=(B(a)B(b))\vartriangleright c,
$$
we see that $(H,\c,\circ_B)$ is a $H^{\rm op}$-module Hopf brace if and only if $B(ba)\vartriangleright c=(B(a)B(b))\vartriangleright c$ for all $a,b,c\in H$.

\smallskip

\noindent
\meqref{it:rbshb3}
If $B$ is an anti-algebra homomorphism, then $B(yx)=B(x)B(y)=B(x\circ_B y)$ for all $x,y\in H$.
So for any $a,b,c\in H$, we obtain the desired equality:
\begin{eqnarray*}
	a_{(1)}b_{(1)}((B(a_{(2)}b_{(2)})B(T(a_{(3)})))\vartriangleright c)&=&a_{(1)}b_{(1)}((B(b_{(2)}\circ_B a_{(2)})B(T(a_{(3)})))\vartriangleright c)\\
	&=&a_{(1)}b_{(1)}((B((b_{(2)}\circ_B a_{(2)})\circ_B T(a_{(3)})))\vartriangleright c)\\
	&=&ab_{(1)}(B(b_{(2)})\vartriangleright c).
\end{eqnarray*}
Then by \meqref{it:rbshb1}, $(H, \c, \circ_B)$ is a symmetric Hopf brace.
\end{proof}

\noindent {\bf Acknowledgments. }  This work was supported by National Natural Science Foundation of China (12201188), Fundamental Research Funds for the Central Universities (ZJ22195010), and postdoctoral research grant in Henan Province (No.202103090).

\noindent
{\bf Declaration of interests. } The authors have no conflicts of interest to disclose.

\noindent
{\bf Data availability. } Data sharing is not applicable as no new data were created or analyzed.

\end{document}